\documentclass[12pt,a4paper,oneside]{amsart}
\usepackage{amsfonts, amsmath, amssymb, amsthm, amscd}
\usepackage{hyperref}
\hypersetup{
  colorlinks   = true, 
  urlcolor     = blue, 
  linkcolor    = blue, 
  citecolor   = red 
}
\usepackage{anysize}
\marginsize{2cm}{2cm}{1.5cm}{1.5cm}

\usepackage[T2A]{fontenc}
\usepackage[cp1251]{inputenc}
\usepackage[english]{babel}
\usepackage{cmap}
\usepackage{enumerate}
\usepackage[pdftex]{graphicx}

\newtheorem{theorem}{Theorem}[section]
\newtheorem{lemma}[theorem]{Lemma}

\newtheorem{corollary}[theorem]{Corollary}

\theoremstyle{definition}

\theoremstyle{remark}
\newtheorem{remark}[theorem]{Remark}

\numberwithin{equation}{section}

\newcommand{\R}{\mathbb{R}}

\def\Const{2e^2}%

\newcommand{\mgl}[2]{\rho_{#1}\!\left( #2 \right)}
\newcommand{\mglx}[1]{\rho_{X}\!\!\:\left( #1 \right)}%
\newcommand{\rmglx}[1]{\rho^{-1}_{X}\!\!\!\:\left( #1 \right)}%
\newcommand{\dist}[2]{{\rm dist}\left(#1, #2 \right)}

\newcommand{\norm}[1]{\left\| #1 \right\|}

\def\co{\mathop{\rm conv}}

\def\diam{\mathop{\rm diam}}

\newcommand{\iprod}[2]{\left\langle#1,#2\right\rangle}

\newcommand{\prooff}{\noindent {\bf Proof.}\\}%
\newcommand{\bbox}{\par\noindent\ensuremath{\Box}\par\noindent}%
\newenvironment{prf}
{\prooff}
{\bbox}

\title{Approximate Carath{\'e}odory's theorem in uniformly smooth Banach spaces}

\author{Grigory Ivanov}
\address{Department of Mathematics\\
University of Fribourg\\
Chemin du Mus\'ee 23\\
CH-1700 Fribourg P\'erolles\\
SWITZERLAND}
\address{Department of Higher Mathematics, Moscow Institute of Physics and Technology,  Institutskii pereulok 9, Dolgoprudny, Moscow
region, 141700, Russia.}
\email{grimivanov@gmail.com}
\thanks{Supported by the Swiss National Science Foundation grant 200021\_179133. Supported by Russian Foundation for Basic Research, project 18-01-00036 A. Research partially supported by Swiss National Science Foundation grants 200020-165977 and 200021-162884}

\begin{document}
\maketitle

{\bf Abstract.} We study the 'no-dimension' analogue of Carath{\'e}odory's theorem in  Banach spaces. We prove  such a result together with its colorful version for uniformly smooth Banach spaces. It follows that uniform smoothness leads to a greedy de-randomization of   Maurey's classical lemma \cite{pisier1980remarques}, which  is itself a 'no-dimension' analogue of Carath{\'e}odory's theorem with a probabilistic proof. We find the asymptotically tight upper bound on the deviation of the convex hull  from the $k$-convex hull of a bounded set in $L_p$ with $1 < p \leq 2$ and get asymptotically the same bound as in Maurey's lemma  for $L_p$  with $1 < p < \infty.$  \nocite{cassels1975measures}
\bigskip


\section{Introduction}  
Carath{\'e}odory's theorem \cite{caratheodory1907variabilitatsbereich} is a classical result in Convex Geometry. It states that every point in  the convex hull $\co S$ of a set $S \subset \R^d$  is a convex combination of at most $d+1$ points of $S.$ There are many different generalizations of this theorem. Recently, several authors studied so-called approximate versions of Carath{\'e}odory's theorem, where the distance between a point of the convex hull $\co S $ of a bounded set $S$ and the $k$-convex hull $\co_k S$ were investigated. The latter, the $k$-convex hull of $S,$ is the set of all convex combinations of at most $k$ points of $S.$

For example, in Theorem 2.2 of  \cite{adiprasito2019theorems} the following optimal result in the Euclidean case was proven:

{\it The distance between  any point $p$ in the convex hull of a bounded set $S$ of a Euclidean space  and  the $k$-convex hull $\co_k S$ is at most $\frac{\diam S}{\sqrt{2k}}$.} \\ \noindent
However, not only the Euclidean case was studied. 
Probably, the most significant result in the area is Maurey's lemma \cite{pisier1980remarques}, in which   an approximate version of Carath{\'e}odory's theorem is proven   for spaces that have (Rademacher) type $p > 1$. We explain the definitions and formulate Maurey's lemma and explain its relation to our results in the next  section. We note here that Maurey's lemma is a more general statement than our results.  
However, Maurey's proof uses Khintchine's inequality and it  is probabilistic. 

We think that a simple geometric property of  Banach spaces hides behind   such a sophisticated technique.    
In this paper, we prove the approximate Carath{\'e}odory's theorem for uniformly smooth Banach spaces, and bound the distance between a point of the convex hull  of a bounded set and its $k$-convex hull in terms of the modulus of smoothness of a Banach space. In fact, we provide a greedy algorithm for the approximation of a point  of the convex hull of a bounded set in a uniformly smooth Banach space. We note here that an $L_p$ space for  $1 < p < \infty$ is uniformly smooth and its modulus of smoothness is well-known.  We give all the required definitions in Section \ref{sec:mgl} and briefly discuss there some  properties of Banach spaces that are connected to the uniform smoothness of Banach space.
The following statement is the main result of the paper.
\begin{theorem}
\label{thm:no-dim_caratheodory}
Let $S$ be a bounded set in a uniformly smooth Banach space $X,$ $a \in \co S.$
Then there exists a sequence $\{x_i\}_1^\infty \subset S$ such that for  vectors
$a_k = \frac{1}{k} \sum\limits_{i \in [k]} x_i$ the following inequality holds
\begin{equation} \label{eqn:ineq_main_nodim_caratheodory_or}
\norm{a - a_k} \leq \frac{\Const}{k \rmglx{1/k}} \diam S,
\end{equation}
where $\mglx{\cdot}$ is the modulus of smoothness of $X.$
\end{theorem}
In the proof we provide a  greedy algorithm for constructing such a sequence,
which directly implies the colorful version of approximate Carath{\'e}odory's theorem.
\begin{corollary}\label{cor:colorful_caratheodory}
Let $\{S_i\}_1^\infty$ be a family of sets of a Banach space $X$ such that 
$a \in \bigcap\limits_1^{\infty} \co S_i$ and $D = \sup_i \diam S_i < \infty.$
Then there exists a transversal sequence $\{x_i\}_1^\infty$ (that is, $x_i \in S_i$) such that for  vectors
$a_k = \frac{1}{k} \sum\limits_{i \in [k]} x_i$ the following inequality holds
\begin{equation}
\label{eqn:ineq_main_nodim_caratheodory}
\norm{a - a_k} \leq \frac{\Const}{k \rmglx{1/k}} D,
\end{equation}
where $\mglx{\cdot}$ is the modulus of smoothness of $X.$
\end{corollary}

At the end of the next section, we show that the bound in inequality \eqref{eqn:ineq_main_nodim_caratheodory_or} is  tight  for all 
$L_p$ spaces with $ 1 \leq p \leq 2$ and coincides with the bound obtained in \cite{barman2015approximating} for $2 \leq p < \infty$ up to a reasonable constant factor.

There are some others "measures of non-convexity" of the $k$-convex hull. We refer the reader to the survey \cite{fradelizi2017convexification}. However, they were mostly studied in the Euclidean case.

\section{Modulus of smoothness and its properties}
\label{sec:mgl}
In this section, we give definitions from the Banach space theory and give a certain reformulation of the main result using this language.

The modulus of smoothness (or Lindenstrauss' modulus) of a Banach space $X$ is the function $\rho_X: [0, \infty] \to [0, \infty]$ given by
$$
	\mglx{\tau}= 
	\sup \left\{\frac{1}{2}(\norm{x + \tau y} + \norm{x - \tau y}) -1 \quad 
	\middle\vert\
	\quad \norm{x} = \norm{y} \leq  1  \right\}.
$$ 
A Banach space $X$ is called {\it uniformly smooth} if 
$\mglx{t} = o(t)$ as $t \to 0.$ 
It is known that uniform smoothness is equivalent to the uniform differentiability of the norm.
As a good reference with simple geometric explanations of different properties of 
the modulus of smoothness and uniformly smooth spaces we refer to Chapter 2 of \cite{DiestelEng}.

It is known that the modulus of smoothness $\mglx{\cdot}$ of a Banach space  $X$ is a convex strictly increasing function that satisfies the following inequality of  Day-Nordlander type \cite{lindenstrauss1963modulus} for all positive $\tau$ 
\[
\sqrt{1 + \tau^2} -1 = \mgl{H}{\tau} \leq \mglx{\tau},
\]
where $\mgl{H}{\tau}$ is the modulus of smoothness of Hilbert space.
The latter inequality implies the following technical observations, which we use in the sequel,
\begin{equation}
\label{eqn:ineq_smt_1_2}
\frac{\Const}{\rmglx{1}} \geq 1; \quad \frac{\Const}{\rmglx{1/2}} \geq 2;
\end{equation}
and
\begin{equation}
\label{eqn:ineq_smt_k}
\frac{1}{\rmglx{1/k}} \geq 1  \quad \text{for} \quad k \geq 3.
\end{equation}

Clearly, uniform smoothness of the norm is not stable under small perturbations of the norm. However, equivalent renormalization of a space just adds a constant factor to  inequalities \eqref{eqn:ineq_main_nodim_caratheodory_or}, \eqref{eqn:ineq_main_nodim_caratheodory}. That is, we can approximate a point of the convex hull by a point of the $k$-convex hull of a bounded set in every Banach space, which admits an equivalent uniformly smooth norm. Such spaces are well-studied; moreover, due to Enflo \cite{Enflo1972} and Pisier \cite{pisier1975martingales},
we know that these spaces are exactly the so-called super-reflexive Banach spaces. Summarizing their results,  the following assertions are equivalent
\begin{itemize}
\item a Banach space $X$ admits an equivalent uniformly smooth norm;
\item a Banach space $X$ admits an equivalent uniformly smooth norm with modulus of smoothness of {\it power type} $p;$
\item a Banach space $X$  is a super-reflexive space. 
\end{itemize}
It is said that a uniformly smooth Banach space $X$ has modulus of smoothness of
{\it power type} $p$ if, for some $0 < C < \infty,$ $\mglx{\tau} \leq C \tau^p.$

Now we can reformulate our results in a norm renormalization invariant way.
We use $\dist{a}{B}$ to denote the distance between a point $a$ and a set $B.$  
We define the deviation of a set $A$ from a  set $B$ as follows
\[
h^+ (A, B) = \sup\limits_{a \in A} \dist{a}{B}
\] 
We say that  a Banach space $X$ has  {\it Carath{\'e}odory's approximation property} if there is a function $f: \mathbb{N} \to [0, \infty)$ with  
$\lim\limits_{k \to \infty} f(k) = 0$ such that 
$h^{+} (\co S, \co_k S) \leq f(k)$ for  an arbitrary subset $S$ of the unit ball of $X.$ 

Theorem \ref{thm:no-dim_caratheodory} and the discussion above imply the following. 
\begin{corollary}\label{cor:rate}
A super-reflexive Banach  space $X$ has  Carath{\'e}odory's approximation property: that is, there exists $p \in (1,2]$   and   constant $C$ such that
$$
	h^{+} (\co S, \co\nolimits_k S) \leq \varepsilon,
$$  
for all $\varepsilon > 0$ and $k \geq C \left(\frac{\diam S}{\varepsilon}\right)^{\frac{p}{p-1}}.$
\end{corollary}

A Banach space $X$ is said to be of {\it type} $p$ for some $1 < p \leqslant 2,$ if there exists a constant $T_p(X) < \infty$ so that, for every finite set of vectors $\left\{x_{j}\right\}_{j=1}^{n}$ 
in $X,$ we have 
$$
\int_{0}^{1}\left\|\sum_{j=1}^{n} r_{j}(t) x_{j}\right\| d t \leqslant T_p(X)\left(\sum_{j=1}^{n}\left\|x_{j}\right\|^{p}\right)^{1 / p},
$$
where  $\left\{r_{j}\right\}_{j=1}^{\infty}$ denotes the sequence of the Rademacher functions.

We can formulate Maurey's lemma as follows (see also Lemma D in \cite{bourgain1989duality}).

{\it 
Let $S$ be a bounded set in a  Banach space $X$ which is of type $p,$ $a \in \co S.$
Set $q = \frac{p}{p -1}.$
Then there exists a sequence $\{x_i\}_1^\infty \subset S$ such that for  vectors
$a_k = \frac{1}{k} \sum\limits_{i \in [k]} x_i$ the following inequality holds
\begin{equation*} 
\norm{a - a_k} \leq T_p (X) k^{- 1/r} \diam S.
\end{equation*}
}

Since a Banach space with modulus of smoothness of power type p implies type p, then 
Theorem \ref{thm:no-dim_caratheodory} follows from Maurey's lemma (up to a constant term in the inequalities). The converse is not true in general (see \cite{james1978nonreflexive}, \cite{pisier1987random}).  However, type $p$ with some additional not very restrictive property (see table on p.~101 in \cite{lindenstrauss2013classical}) implies modulus of smoothness of power type $p.$

It is known  \cite{lindenstrauss1963modulus} that $L_p,$ $1 < p < \infty,$ is of type $q = \max\{p, 2\}$ and  has modulus of smoothness of power $\max\{p, 2\}.$ More precisely,  
\begin{equation}
\label{mod_smooth_p2infty}
\mgl{L_p} {t} = \left(\frac{(1 +t)^p + (1 -t)^p)}{2}\right)^{1/p} -1 < 
 \frac{p-1}{2} t^2  \quad \text{for} \quad   2 \leqslant p < \infty;
\end{equation}
\begin{equation}
\label{mod_smooth_p12}
\mgl{L_p} {t} = (1 +t^p)^{1/p} -1 \leq \frac{t^p}{p} \quad \text{for} \quad 1 < p \leqslant 2.
\end{equation}
Since $p$ in Corollary \ref{cor:rate} can be chosen to be the order of the modulus of smoothness, we see that $k = O \left(\frac{\diam^2 S}{\varepsilon^2} \right)$ in $L_p,$ $2 \leq p < \infty,$ which coincides  with the rate of convergence in Maurey's lemma. Moreover, as was shown in Theorem 3.2 of Barman's paper \cite{barman2015approximating}, the constants in both inequalities are close enough to one another in this case.  

In case of an $L_p$ space with $1 \leq p  \leq 2,$ the bound in inequality \eqref{eqn:ineq_main_nodim_caratheodory_or} is the best possible up to a constant:
when $n = 2k,$ $S=\{e_1, \dots, e_{2n}\}$ and $a = \{1/n, \dots, 1/n\},$
then $\dist{a}{\co_k S} \geq \frac{1}{4} k^{1/p - 1}.$ Moreover, this implies that $L_1$ and $L_\infty$ have no Carath{\'e}odory's approximation property.

\section{Geometrical idea behind the proof}
We begin with a sketch of the proof of the following folklore analogue of the main theorem in a Euclidean space (see  \cite{starr1969quasi}, \cite{cassels1975measures}).

{\it 
Let $S$ be a subset of  the unit ball of a Euclidean space  such that $0 \in \co S.$
Then there exists a sequence $\{x_i\}_1^\infty \subset S$ such that 
\[ 
\norm{\sum\limits_1^k x_i} \leqslant \sqrt{k}
\]
for all  $k \in \mathbb{N}.$
}  

We apply induction on $k.$ 
Assume that we have chosen $x_1, \dots, x_n$ that satisfy the inequality for all $k \in [n].$ 
Since $0 \in \co S,$ there exists $x_{n+1}$ such that  $\iprod{x_{n+1}}{\sum\limits_1^n x_i}  \leq 0.$  Then, by the law of cosines,
\[
	\norm{\sum\limits_1^n x_i + x_{n+1}} \leq \sqrt{\norm{\sum\limits_1^n x_i}^2+ \norm{x_{n+1}}^2} \leq \sqrt{n+1}.
\]
The statement is proven.

In fact, the law of cosines can be reformulated in terms of the deviation of the unit sphere from its supporting hyperplane   as follows.

We use $u^\ast$ to denote a unit functional that attains its norm on a non-zero vector $u$ of a Banach space $X,$ i.e. $\iprod{u^\ast}{u} = \norm{u} \norm{u^\ast} = \norm{u}.$ Clearly, a set $\{x \in X \mid \iprod{u^\ast}{x} = 1 \}$ is a supporting hyperplane to the unit ball of $X$ at $u/\norm{u}.$    For simplicity, we assume that $u^\ast = 0$ for $u =0.$

Let $H$ be a supporting hyperplane at a unit vector $u$ of the unit ball in a Euclidean space $E$ and $x$ be such that $\iprod{u^\ast}{x} \leq 0.$
Then the norm of $\norm{u+x}$ is at most 
$ \sqrt{1 +\norm{x}^2}$ ($\approx 1 + \mgl{E}{\norm{x}}$ for sufficiently small $\norm{x}$)  with equality only for $x$ parallel to $H.$
That is,  we  see that  vectors from the supporting hyperplane spoil the sum most badly and we measure the deviation of the unit sphere from a supporting hyperplane  to bound the norm  of the sum on each step.
Similar statements can be proven in a uniformly smooth Banach space. However, 
it is not necessarily true that for an arbitrary unit vector $u$  of a Banach space the maximum of norm $\norm{u + v},$ where $\norm{v}$ is fixed and $\iprod{u^\ast}{v} \leq 0,$  is attained on the supporting hyperplane to the unit ball at $u$. 
But one can measure this deviation using the modulus of smoothness, which we do in the following simple lemma.

\begin{lemma}
\label{lem:geom_lemma}
Let $u$ be a unit vector of Banach space $X$ and $x$ such that $\iprod{u^\ast}{x} \leq 0.$  Then 
\[
	\norm{u + x} \leq 2 \mglx{\norm{x}} + 1.
\]
\end{lemma}
\begin{prf}
Since $\iprod{u^\ast}{x} \leq 0$ and $H_u = \{q \in X \mid \quad  \iprod{u^\ast}{q} = 1\}$ is a supporting hyperplane for the unit ball of $X,$ we have that $\norm{u - x} \geq 1.$ Therefore, by the definition of the modulus of smoothness, we obtain
\[
\norm{u +x} + \norm{u-x} \leq 2 \mglx{\norm{x}} + 2,
\]
or, equivalently,
\[
\norm{u +x}  \leq 2 \mglx{\norm{x}} + 2 -  \norm{u-x} \leq  2 \mglx{\norm{x}} + 1.
\]
\end{prf}

Our proof follows the same line as in the above-mentioned Euclidean analogue with the use of Lemma \ref{lem:geom_lemma} instead of the law of cosines.
\begin{remark}
The deviation of the unit sphere from the supporting hyperplane in a Banach space was studied by the author in \cite{Ivanov_supp_modulus}. In particular, it was shown that the inequality from Lemma \ref{lem:geom_lemma} is asymptotically tight as $\norm{x}$ tends to zero. On the other hand, it was proven in \cite{ivanov2017new} that for any $\tau$  in an arbitrary Banach space,  there is a unit vector $u$ and a vector $x$  parallel to a supporting hyperplane  to the unit ball at $u$ such that $\norm{u + x} = \sqrt{1 + \norm{x}^2} = \sqrt{1 + \tau^2}.$ That is, the smallest upper bound is attained in the  Euclidean case.  
\end{remark}

\section{Proof of the main result}
By setting $S_i = S$ for all $i \in \mathbb{N}$, we see that  Corollary \ref{cor:colorful_caratheodory} implies Theorem \ref{thm:no-dim_caratheodory}.
We give the proof of the  Corollary, which  coincides with the proof of the Theorem up to above-mentioned  renaming of sets. \\ \noindent
\begin{prf}
To simplify the proof, we translate $S_i$ to $S_i - a$ and then scale them in such a way that $D = \sup_i\diam S_i = 1.$ Then inequality \eqref{eqn:ineq_main_nodim_caratheodory} transforms into
\begin{equation}
\label{eqn:ineq_main_nodim_caratheodory_simplification}
	\norm{a_k} \leq \frac{\Const}{k \rmglx{1/k}}.
\end{equation}
Let us construct a sequence $\{x_i\}_1^\infty, $ where $x_i \in S_i,$ that satisfies the latter inequality. 
We use the following {\bf algorithm}:
\begin{enumerate}
\item $x_1$ is an arbitrary point from $S_1 .$ 
\item for the  constructed sequence $\{x_1, \dots, x_{k}\},$ $k \geq 2,$ we choose $x_{k+1} \in S_{k+1}$ such that 
$\iprod{x_{k+1}}{a_k^\ast} \leq 0$
(that is, $x_{k+1}$ is an arbitrary point of $S_{k+1}$ if $a_k = 0$).
\end{enumerate}

Firstly, the sequence $\{x_i\}_1^\infty$ is well-defined. 
Indeed,  there exists $q \in S_{k+1}$ such that $\iprod{q}{p} \leq 0$ for an arbitrary functional $p \in X^\ast,$ since $0 \in \co S_{k+1}.$
In the algorithm we choose a
functional $p = a_k^\ast$
such that $\norm{a_k^\ast} = 1$ and  for $a_k = \frac{1}{k}\sum\limits_{i \in [k]} x_i, 
a_k \neq 0,$ the equality $\iprod{a_k^\ast}{a_k} = \norm{a_k}$ is valid.

Secondly, let us show that $\{x_i\}_1^\infty$ satisfies inequality \eqref{eqn:ineq_main_nodim_caratheodory_simplification}. 
Fix $k.$
By inequality \eqref{eqn:ineq_smt_1_2}, we  assume that $k \geq 3.$
By definition put 
$\eta_k =  1 /\rmglx{1/k}$ and $u_k = k a_k.$  
If  inequality $\norm{u_k} \leq \eta_k$ holds, it implies the required inequality. Assume that $\norm{u_k} > \eta_k.$ Then let $m-1$ be the biggest numbers in $[k-1]$ such that $\norm{u_{m-1}} \leq \eta_k$ and $\norm{u_m} > \eta_k.$ Such $m$ exists since, by inequality \eqref{eqn:ineq_smt_k}, we have that $\norm{u_1} = \norm{x_1} \leq \eta_k.$  Since $\eta_k  \geq 1$ for $k \geq 3$, we have 
\begin{equation}
\label{eqn:u_m}
	\norm{u_m} \leq \eta_k + 1 \leq 2 \eta_k.
\end{equation}
By the choice of $m,$ we have that 
$\norm{u_{m+\ell-1}} > \eta_k \geq 1$
for all $\ell \in [k-m].$ 
Hence, by the definition of $\eta_k$ and the construction of $\{x_i\}_1^\infty,$ we get
\begin{gather*}
	\norm{u_{m+\ell}} = \norm{u_{m + \ell - 1} + x_{m + \ell}} = 
	\norm{u_{m + \ell - 1}} \norm{\frac{u_{m + \ell - 1}}{\norm{u_{m + \ell - 1}}} + \frac{x_{m + \ell}}{\norm{u_{m + \ell - 1}}}} 
	\stackrel{\text{Lemma} \; \ref{lem:geom_lemma}}{\leq} \\
	\norm{u_{m + \ell - 1}} \left(1 + 2 \mglx{\frac{\norm{x_{m + \ell}}}{\norm{u_{m + \ell - 1}}}} \right) \leq 
		\norm{u_{m + \ell - 1}} \left(1 + 2 \mglx{\frac{1}{\eta_k}} \right) = 
		\norm{u_{m + \ell - 1}} \left(1 + \frac{2}{k} \right).
\end{gather*}
Combining these inequalities for all $\ell \in [k-m]$  and inequality \eqref{eqn:u_m}, we obtain
\[
\norm{u_k} \leq \norm{u_m} \left(1 + \frac{2}{k} \right)^{k-m} \leq \Const \eta_k.
\]
Therefore, $a_k = u_k/k$ satisfies inequality \eqref{eqn:ineq_main_nodim_caratheodory_simplification}.
The theorem is proven.  
\end{prf}

We note here that the authors used a similar idea in \cite{ivanov2012generalization} to prove the convexity of a special type limit object in a uniformly smooth Banach space.

{\bf Acknowledgements.}  
I wish to thank Imre B{\'a}r{\'a}ny for bringing the problem to my attention. 
\bibliographystyle{abbrv}
\bibliography{/home/grigory/Dropbox/work_current/uvolit} 

\begin{thebibliography}{10}

\bibitem{adiprasito2019theorems}
K.~Adiprasito, I.~B{\'a}r{\'a}ny, and N.~H. Mustafa.
\newblock Theorems of {C}arath{\'e}odory, {H}elly, and {T}verberg without
  dimension.
\newblock In {\em Proceedings of the Thirtieth Annual ACM-SIAM Symposium on
  Discrete Algorithms}, pages 2350--2360. SIAM, 2019.

\bibitem{barman2015approximating}
S.~Barman.
\newblock Approximating {N}ash equilibria and dense bipartite subgraphs via an
  approximate version of {C}aratheodory's theorem.
\newblock In {\em Proceedings of the forty-seventh annual ACM symposium on
  Theory of computing}, pages 361--369. ACM, 2015.

\bibitem{bourgain1989duality}
J.~Bourgain, A.~Pajor, S.~J. Szarek, and N.~Tomczak-Jaegermann.
\newblock On the duality problem for entropy numbers of operators.
\newblock In {\em Geometric aspects of functional analysis}, pages 50--63.
  Springer, 1989.

\bibitem{caratheodory1907variabilitatsbereich}
C.~Carath{\'e}odory.
\newblock {\"U}ber den variabilit{\"a}tsbereich der koeffizienten von
  potenzreihen, die gegebene werte nicht annehmen.
\newblock {\em Mathematische Annalen}, 64(1):95--115, 1907.

\bibitem{cassels1975measures}
J.~Cassels.
\newblock Measures of the non-convexity of sets and the
  {S}hapley--{F}olkman--{S}tarr theorem.
\newblock In {\em Mathematical Proceedings of the Cambridge Philosophical
  Society}, volume~78, pages 433--436. Cambridge University Press, 1975.

\bibitem{DiestelEng}
J.~Diestel.
\newblock {\em Geometry of Banach Spaces - Selected Topics}, volume 485.
\newblock Springer-Verlag Berlin Heidelberg, 1975.

\bibitem{Enflo1972}
P.~Enflo.
\newblock Banach spaces which can be given an equivalent uniformly convex norm.
\newblock {\em Isr. J. Math.}, 13:281--288, 1972.

\bibitem{fradelizi2017convexification}
M.~Fradelizi, M.~Madiman, A.~Marsiglietti, and A.~Zvavitch.
\newblock The convexification effect of minkowski summation.
\newblock {\em arXiv preprint arXiv:1704.05486}, 2017.

\bibitem{Ivanov_supp_modulus}
G.~M. Ivanov.
\newblock Modulus of supporting convexity and supporting smoothness.
\newblock {\em Eurasian Math. J.}, 6(1):26--40, 2015.

\bibitem{ivanov2017new}
G.~M. Ivanov and H.~Martini.
\newblock New moduli for {B}anach spaces.
\newblock {\em Annals of Functional Analysis}, 8(3):350--365, 2017.

\bibitem{ivanov2012generalization}
G.~M. Ivanov and E.~S. Polovinkin.
\newblock A generalization of the set averaging theorem.
\newblock {\em Mathematical Notes}, 92(3-4):369--374, 2012.

\bibitem{james1978nonreflexive}
R.~James.
\newblock Nonreflexive spaces of type 2.
\newblock {\em Israel Journal of Mathematics}, 30(1-2):1--13, 1978.

\bibitem{lindenstrauss1963modulus}
J.~Lindenstrauss.
\newblock On the modulus of smoothness and divergent series in {B}anach spaces.
\newblock {\em The Michigan Mathematical Journal}, 10(3):241--252, 1963.

\bibitem{lindenstrauss2013classical}
J.~Lindenstrauss and L.~Tzafriri.
\newblock {\em Classical Banach spaces II: function spaces}, volume~97.
\newblock Springer Science \& Business Media, 2013.

\bibitem{pisier1975martingales}
G.~Pisier.
\newblock Martingales with values in uniformly convex spaces.
\newblock {\em Israel Journal of Mathematics}, 20(3-4):326--350, 1975.

\bibitem{pisier1980remarques}
G.~Pisier.
\newblock Remarques sur un r{\'e}sultat non publi{\'e} de {B}. {M}aurey.
\newblock {\em S{\'e}minaire Analyse fonctionnelle (dit" Maurey-Schwartz")},
  pages 1--12, 1980.

\bibitem{pisier1987random}
G.~Pisier and Q.~Xu.
\newblock Random series in the real interpolation spaces between the spaces
  $v_p$.
\newblock In {\em Geometrical Aspects of Functional Analysis}, pages 185--209.
  Springer, 1987.

\bibitem{starr1969quasi}
R.~M. Starr.
\newblock Quasi-equilibria in markets with non-convex preferences.
\newblock {\em Econometrica: journal of the Econometric Society}, pages 25--38,
  1969.

\end{thebibliography}
\end{document}